\def\version{0.59}

\def\journal{JA or JKT}
\def\titlep{Some inverse limits of Cuntz algebras}
\documentclass[11pt,epic,eepic]{article}
\usepackage{graphicx,ifthen}
\usepackage{amssymb}
\usepackage{amsmath}
\usepackage{color}
\usepackage{epic,eepic}
\usepackage{latexsym}

 at12pt

\newcommand{\qed}{\hbox{\rule[-2pt]{3pt}{6pt}}}
\newcommand{\qedh}{\hfill\qed \\}

\newcommand{\vv}{\vspace{.3in}}

\newtheorem{Thm}{Theorem}[section]
\newtheorem{que}[Thm]{Question}

\newtheorem{rem}[Thm]{Remark}

\newtheorem{defi}[Thm]{Definition}
\newtheorem{lem}[Thm]{Lemma}

\newtheorem{prop}[Thm]{Proposition}

\newtheorem{cor}[Thm]{Corollary}
\newtheorem{fact}[Thm]{Fact}
\newtheorem{fig}[Thm]{Figure}

\newcommand{\ww}{\vv\noindent}

\newcommand{\kn}{\Large\bf
$K\hspace{-.4cm} N$
\Large\bf\vv }

\def\cal#1{\mathcal #1}
\def\con{{\cal O}_{n}}
\def\coni{{\cal O}_{\infty}}

\def\pr{{\it Proof.}\quad}

\def\co#1{{\cal O}_{#1}}
\def\disp#1{{\displaystyle #1}}
\def\brl{branching law}
\def\bfsnl{{\rm BFS}_{N}(\Lambda)}

\addtocounter{footnote}{1}
\def\sftt#1{
\setcounter{equation}{0}
\addtocounter{footnote}{1}
\section{#1}
}

\def\ssft#1{\subsection{#1}}
\def\sssft#1{\subsubsection{#1}}
\def\ssfr#1{\subsection*{#1}}

\begin{document}
%
%
\def\autherp{Katsunori Kawamura}
\def\emailp{e-mail: kawamura@kurims.kyoto-u.ac.jp.}
\def\addressp{{\small {\it College of Science and Engineering,
 Ritsumeikan University,}}\\
{\small {\it 1-1-1 Noji Higashi, Kusatsu, Shiga 525-8577, Japan}}
}

\def\infw{\Lambda^{\frac{\infty}{2}}V}
\def\zhalfs{{\bf Z}+\frac{1}{2}}
\def\ems{\emptyset}
\def\pmvac{|{\rm vac}\!\!>\!\! _{\pm}}
\def\vac{|{\rm vac}\rangle _{+}}
\def\dvac{|{\rm vac}\rangle _{-}}
\def\ovac{|0\rangle}
\def\tovac{|\tilde{0}\rangle}
\def\expt#1{\langle #1\rangle}
\def\zph{{\bf Z}_{+/2}}
\def\zmh{{\bf Z}_{-/2}}
\def\brl{branching law}
\def\bfsnl{{\rm BFS}_{N}(\Lambda)}
\def\scm#1{S({\bf C}^{N})^{\otimes #1}}
\def\mqb{\{(M_{i},q_{i},B_{i})\}_{i=1}^{N}}
\def\zhalf{\mbox{${\bf Z}+\frac{1}{2}$}}
\def\zmha{\mbox{${\bf Z}_{\leq 0}-\frac{1}{2}$}}
\newcommand{\mline}{\noindent
\thicklines
\setlength{\unitlength}{.1mm}
\begin{picture}(1000,5)
\put(0,0){\line(1,0){1250}}
\end{picture}
\par
 }
\def\ptimes{\otimes_{\varphi}}
\def\qtimes{\otimes_{\tilde{\varphi}}}
\def\delp{\Delta_{\varphi}}
\def\delps{\Delta_{\varphi^{*}}}
\def\gamp{\Gamma_{\varphi}}
\def\gamps{\Gamma_{\varphi^{*}}}
\def\sem{{\sf M}}
\def\sen{{\sf N}}
\def\hdelp{\hat{\Delta}_{\varphi}}
\def\tilco#1{\tilde{\co{#1}}}
\def\tilconi{\tilde{{\cal O}}_{\infty}}
\def\hatconi{\hat{{\cal O}}_{\infty}}
\def\hatcola{\hat{{\cal O}}(\Lambda)}
\def\ba{\mbox{\boldmath$a$}}
\def\bb{\mbox{\boldmath$b$}}
\def\bc{\mbox{\boldmath$c$}}
\def\be{\mbox{\boldmath$e$}}
\def\bp{\mbox{\boldmath$p$}}
\def\bq{\mbox{\boldmath$q$}}
\def\bu{\mbox{\boldmath$u$}}
\def\bv{\mbox{\boldmath$v$}}
\def\bw{\mbox{\boldmath$w$}}
\def\bx{\mbox{\boldmath$x$}}
\def\by{\mbox{\boldmath$y$}}
\def\bz{\mbox{\boldmath$z$}}
\def\titlepage{

\noindent
{\bf 
\noindent
\thicklines
\setlength{\unitlength}{.1mm}
\begin{picture}(1000,0)(0,-300)
\put(0,0){\kn \knn\, for \journal\, Ver.\version}
\put(0,-50){\today}
\end{picture}
}
\vspace{-2.3cm}
\quad\\
{\small file: \textsf{tit01.txt,\, J1.tex}
\footnote{
${\displaystyle
\mbox{directory: \textsf{\fileplace}, 
file: \textsf{\incfile},\, from \startdate}}$}}
\quad\\
\framebox{
\begin{tabular}{ll}
\textsf{Title:} &
\begin{minipage}[t]{4in}
\titlep
\end{minipage}
\\
\textsf{Author:} &\autherp
\end{tabular}
}
{\footnotesize	
\tableofcontents }
}

%
%
%
\setcounter{section}{0}
\setcounter{footnote}{0}
\setcounter{page}{1}
\pagestyle{plain}

%
%
\title{\titlep}
\author{\autherp\thanks{\emailp}
\\
\addressp}
\date{}
\maketitle
%
%
\begin{abstract}
We construct a nontrivial inverse system of 
Cuntz algebras $\{{\cal O}_{n}:2\leq n<\infty\}$, whose
inverse limit is $*$-isomorphic onto ${\cal O}_{\infty}$.
By using this result,
it is shown that
the $K_{0}$-functor is discontinuous with respect to
the inverse limit
even if the limit is a C$^{*}$-algebra.
\end{abstract}


%
%
\sftt{Introduction}
\label{section:first}
In this paper,
we construct a nontrivial inverse system of 
Cuntz algebras $\{{\cal O}_{n}:2\leq n<\infty\}$, whose
inverse limit is ${\cal O}_{\infty}$:
%
%
\begin{equation}
\label{eqn:firstlimit}
\varprojlim {\cal O}_{n}\cong {\cal O}_{\infty}.
\end{equation}
In order to explain (\ref{eqn:firstlimit}),
we will recall inverse system of C$^{*}$-algebras and 
pro-C$^{*}$-algebra, and 
show the construction of the inverse system in this section.
%
%
\ssft{Unital $*$-homomorphisms among Cuntz algebras}
\label{subsection:firstone}
Our study is motivated by well-known 
facts of unital $*$-homomorphisms among Cuntz algebras.
Hence we start with their explanation in this subsection.
For unital C$^{*}$-algebras $A$ and $B$,
let ${\rm Hom}(A,B)$ denote the set of all unital $*$-homomorphisms
from $A$ to $B$ and let $K_{0}(A)$ denote the $K_{0}$-group of $A$
\cite{BlackadarK}.
%
%
\begin{fact}
\label{fact:ktheorys}
For $2\leq n\leq \infty$,
let $\con$ denote the Cuntz algebra \cite{Cuntz}.
Then the following holds:
\begin{enumerate}
\item
For $2\leq n\leq \infty$ and a unital C$^{*}$-algebra $A$,
if $f\in {\rm Hom}(A,\co{n})$,
then the induced homomorphism $\hat{f}$ from $K_{0}(A)$
to $K_{0}(\con)$
is surjective.
\item
For $2\leq m,n<\infty$, ${\rm Hom}(\co{m},\con)\ne \emptyset$ if and only if
there exists a positive integer $k$ such that $m=(n-1)k+1$.
\item
For any $2\leq m<\infty$, ${\rm Hom}(\co{m},\coni)=\emptyset$.
\item
For any $2\leq n\leq \infty$, ${\rm Hom}(\coni,\con)\ne \emptyset$.
\end{enumerate}
\end{fact}
%
%
\pr
(i)
From Example 6.3.2 in \cite{BlackadarK}, 
the class of the unit of $\con$ is the generator of 
$K_{0}(\con)$.
Since $\hat{f}$ maps the class of the unit of $A$ to
that of $\con$,
the statement holds. 

\noindent
(ii)
This holds from  Lemma 2.1 in \cite{SE01}
(see also \cite{Davidson}, p164,V.16.) 
In $\S$ \ref{subsection:firstthree},
we will given
concrete $*$-homomorphisms among Cuntz algebras.

\noindent
(iii)
From 
\cite{Cuntz1981},
$K_{0}(\con)\cong {\bf Z}/(n-1){\bf Z}$ ($2\leq n<\infty$) and
$K_{0}(\coni)\cong {\bf Z}$.
From these and (i),
the statement holds.

\noindent
(iv)
This will be shown by using concrete $*$-homomorphisms
$\{f_{n,\infty}:n\geq 1\}$ in (\ref{eqn:eone}).
\qedh

\noindent
About homomorphisms among Cuntz algebras,
more general results are known (\cite{GPW}, Lemma 7.1.)

Since $\con$ is simple for each $n$,
any unital $*$-homomorphism from $\con$ is injective,
hence it is an embedding of $\con$. 
For examples of Fact \ref{fact:ktheorys}(ii), the following embeddings 
among $\co{2},\ldots,\co{8}$ (except endomorphisms)
are illustrated as follows (\cite{SE01}, $\S$2.1):

\noindent 
%
\def\setone{
\put(630,870){$\co{2}$}
\put(-200,650){\vector(3,1){500}}
\put(480,650){\vector(1,1){150}}
\put(1150,650){\vector(-2,1){250}}
\put(1700,650){\vector(-3,1){500}}
\put(-400,550){$\co{3}$}
\put(330,550){$\co{4}$}
\put(1170,550){$\co{6}$}
\put(1750,550){$\co{8}$}
\put(-750,450){\vector(3,1){250}}
\put(-1000,350){$\co{5}$}
\put(-80,400){\vector(-2,1){180}}
\put(150,400){\vector(2,1){180}}
\put(-50,300){$\co{7}$}
}
%
%
%
\begin{fig}
\label{fig:three}
\quad\\
\setlength{\unitlength}{.024713mm}
\begin{picture}(4811,800)(-2500,200)
\thicklines
\put(-500,0){\setone}
\end{picture}
\end{fig}

\noindent
where an arrow ``$A\to B$" means a unital embedding of $A$ into $B$.
For $2\leq n<m\leq 8$, there is no unital $*$-homomorphism from $\co{m}$ to $\con$
if there is no oriented path from $\co{m}$ to $\con$ 
in Figure \ref{fig:three}.
%
%
\begin{rem}
\label{rem:zero}
{\rm
In p184 of \cite{Cuntz},
there is a statement about embeddings 
among Cuntz algebras as follows:
%
%
\begin{equation}
\label{eqn:induction}
\co{2}\supset \co{3}\supset \co{4}\supset \ldots \supset \coni.
\end{equation}
It is explained that 
(\ref{eqn:induction}) is given by using the induction for the construction of 
 a certain unital embedding of $\co{3}$ into $\co{2}$.
However, 
(\ref{eqn:induction}) never means unital embeddings because of 
Fact \ref{fact:ktheorys}(ii) except ``$\co{2}\supset \co{3}$"
and ``$\supset \coni$".
}
\end{rem}

On the other hand,  there exist well-known two orders on
the set ${\bf N}$ of all positive integers
$\{1,2,3,\ldots\}$.
The first is the standard linear order $\leq $, that is,
$1\leq 2\leq 3\leq \cdots$.
The second is the order $\preceq$ on ${\bf N}$ (\cite{Bourbaki}, $\S$1.11)
defined  as
%
%
\begin{equation}
\label{eqn:precdef}
m\preceq n  \quad \mbox{ if  $m$ divides $n$. }
\end{equation}
This relation $m\preceq n$ is usually written as $m|n$ in
number theory \cite{HW}.
Both $({\bf N},\leq)$ and $({\bf N},\preceq)$
are directed sets, but
$({\bf N},\preceq)$ is not a
totally ordered set, which 
is illustrated as 
the following directed graph (except relations $n\preceq n$):
%
\def\setone{
\put(630,900){$1$}
\put(300,850){\vector(-3,-1){500}}
\put(600,850){\vector(-1,-1){150}}
\put(800,850){\vector(2,-1){250}}
\put(1100,850){\vector(3,-1){500}}
\put(-400,550){$2$}
\put(330,550){$3$}
\put(1170,550){$5$}
\put(1750,550){$7$}
\put(-500,550){\vector(-3,-1){250}}
\put(-1000,350){$4$}
\put(-260,550){\vector(2,-1){180}}
\put(250,550){\vector(-2,-1){180}}
\put(-50,300){$6$}
}
%
%
\begin{fig}
\label{fig:five}
\quad\\
\setlength{\unitlength}{.024713mm}
\begin{picture}(4811,550)(-2500,350)
\thicklines
\put(-500,0){\setone}
\end{picture}
\end{fig}

\noindent
By comparison with Figure \ref{fig:three},
it is clear that
Figure \ref{fig:five}
is just the graph with inverse direction of 
Figure \ref{fig:three} by rewriting their suffix numbers.
This idea is rigorously verified by using Fact \ref{fact:ktheorys}(ii)
and we can restate Fact \ref{fact:ktheorys}(ii), (iii) and (iv) 
by using the order $\preceq$ as follows.
%
%
\begin{cor}
\label{cor:second}
Let $\hat{{\bf N}}\equiv {\bf N}\cup\{\infty\}$
and extend $\preceq$ on $\hat{{\bf N}}$ as
$\infty \preceq\infty$ and $n\preceq \infty$ for each $n\in {\bf N}$.
Then, 
for $n,m\in \hat{{\bf N}}$,
%
%
\begin{equation}
\label{eqn:homhom}
{\rm Hom}(\co{m+1},\co{n+1})\ne \emptyset\quad
\mbox{if and only if }\quad n\preceq m
\end{equation}
where $n+\infty$ means $\infty$ for convenience.
Especially,
(\ref{eqn:homhom}) holds for each $n,m\in {\bf N}$,
and
there exists no unital $*$-homomorphism
from $\co{m}$ to $\con$ when $n>m$.
\end{cor}
From Corollary \ref{cor:second},
if there exist the following unital inclusions 
%
%
\begin{equation}
\label{eqn:inclusionstwo}
\co{n_{1}+1}\supset \co{n_{2}+1}\supset\co{n_{3}+1}\supset \cdots
\end{equation}
for $\{n_{i}\in {\bf N}:i\geq 1\}$,
then
we obtain order relations
$n_{1}\preceq n_{2}\preceq n_{3}\preceq\cdots$.

It is well-known that 
the order $\preceq$ is used in the theory of profinite groups
\cite{RZ}. 
A {\it profinite group}
is defined as an inverse limit of 
finite groups.
From this and Corollary \ref{cor:second},
the following questions are inspired.
%
%
\begin{que}
\label{que:first}
\begin{enumerate}
\item
Does there exist an inverse system of Cuntz algebras
$\{\co{n+1}:1\leq n<\infty\}$ over the directed set
$({\bf N},\preceq)$ with respect to embeddings in 
Figure \ref{fig:three}?
\item
If such an inverse system in (i) is found,
then what-like algebra is the inverse limit?
\item
If the answer to (ii) is given,
then what is the meaning of this result?
\end{enumerate}
\end{que}

\noindent
The purpose of this paper is to
give answers to these questions.

%
%
\ssft{Inverse limits of C$^{*}$-algebras}
\label{subsection:firsttwo}
In order to consider inverse limits of Cuntz algebras in
Question \ref{que:first},
we recall previous works of 
inverse limits of C$^{*}$-algebras in this subsection.

According to Phillips \cite{Phillips},
Fragoulopoulou \cite{Fragoulopoulou}
and Joi\c{t}a \cite{Joita},
inverse limits (or projective limits 
\cite{Cuntz2005,Palmer1}) of C$^{*}$-algebras
were studied by different names in the literature as follows:
b$^{*}$-algebras (Allan \cite{Allan}, Apostol \cite{Apostol}), 
LMC$^*$-algebras  (Lassner \cite{Lassner}, 
Schmudgen \cite{Schmudgen}), 
l.m.c.C$^*$-algebras
(Mallios \cite{Mallios}),
locally C$^*$-algebras 
(Inoue \cite{Inoue}, 
Fragoulopoulou \cite{Fragoulopoulou1981}),
generalized operator algebras 
(Weidner \cite{Weidner1,Weidner2}),
F$^{*}$-algebras (Brooks \cite{Brooks}),
$\sigma$-C$^{*}$-algebras (Arveson \cite{Arveson}),
or 
pro-C$^{*}$-algebras (Voiculescu \cite{Voiculescu}).

Next,
we recall definitions and basic facts, where
notations are slightly changed from \cite{Phillips}.
Let  $(D,\leq)$ be a {\it directed set} \cite{Kelley}, that is,
$D$ is a non-empty set and $\leq$ is a binary relation on $D$
which satisfies the conditions:
For all $a,b,c\in D$, we have $a\leq a$;
 if $a\leq b$ and $b\leq c$, then $a\leq c$;
if $a,b\in D$, then there exists $c\in D$ such that $a\leq c$ and $b\leq c$.
Such $\leq$ is called a {\it preorder} \cite{Bourbaki}.
We call $\leq$ an order in this paper for simplicity of description.
For every concrete directed set $(D,\leq)$ in this paper,
the order $\leq$ satisfies
the antisymmetric law, that is,
if $a\leq b$ and $b\leq a$, then $a=b$.
%
%
\begin{defi}
\label{defi:inverse}
Let $(D,\leq)$ be a directed set.
\begin{enumerate}
\item
A data $\{(A_{d},\varphi_{d,e}):d,e\in D\}$
is an inverse system (or projective system \cite{Voiculescu}) 
of C$^{*}$-algebras if 
$A_{d}$ is a C$^{*}$-algebra for each $d\in D$ and 
$\varphi_{d,e}$ is a $*$-homomorphism from
$A_{e}$ to $A_{d}$ when $d\leq e$
such that
$\varphi_{d,e}\circ \varphi_{e,f}=\varphi_{d,f}$
when $d\leq e\leq f$,
and $\varphi_{d,d}=id_{A_{d}}$.
\item
The inverse limit $(A,\{\pi_{d}\}_{d\in D})$
of an inverse system 
$\{(A_{d},\varphi_{d,e}):d,e\in D\}$
of C$^{*}$-algebras is 
a topological $*$-algebra $A$
and a $*$-homomorphism $\pi_{d}$
from $A$ to $A_{d}$
such that the following conditions hold:
\begin{enumerate}
\item
$\varphi_{d,e}\circ \pi_{e}=\pi_{d}$ when $d\leq e$,
\item
for any $*$-algebra $B$ and $*$-homomorphisms
$\{\eta_{d}\}_{d\in D}$,
$\eta_{d}:B\to A_{d}$,
 which satisfy $\varphi_{d,e}\circ \eta_{e}=\eta_{d}$ when $d\leq e$,
there exists a unique $*$-homomorphism $\psi$ from $B$ to $A$
such that $\pi_{d}\circ \psi=\eta_{d}$ for each $d\in D$,
\item
the topology of $A$ is the weakest topology 
such that every $\pi_{d}$ is continuous.
\end{enumerate}
In this case,
$A$ is written as 
$\varprojlim_{D}(A_{d},\varphi_{d,e})$ or
$\varprojlim(A_{d},\varphi_{d,e})$,
and $\pi_{d}$ is called the canonical homomorphism (or projection \cite{RZ},
canonical mapping \cite{Bourbaki}.) 
\item
A pro-C$^{*}$-algebra $A$ is a complete Hausdorff topological
$*$-algebra over ${\bf C}$ whose topology is determined by
its continuous C$^{*}$-seminorms
in the sense that a net $\{a_{\lambda}\}$ converges to $0$
if and only if $p(a_{\lambda})\to 0$ for every
continuous C$^{*}$-seminorm $p$ on $A$.
\item
An inverse limit of C$^{*}$-algebras over a countable directed set
is called a $\sigma$-C$^{*}$-algebra.
\end{enumerate}
\end{defi}

\noindent
From ``1.2. Proposition" in \cite{Phillips}, 
a C$^{*}$-algebra $A$ is a pro-C$^{*}$-algebra
if and only if
it is the inverse limit of an inverse system of C$^{*}$-algebras.
Any C$^{*}$-algebra is a pro-C$^{*}$-algebra.
For any inverse system $\{(A_{d},\varphi_{d,e}):d,e\in D\}$,
the inverse limit $\varprojlim(A_{d},\varphi_{d,e})$
is given as the following subset of the product set $\prod_{d\in D}A_{d}$:
%
%
\begin{equation}
\label{eqn:xain}
\{(x_{d})\in \prod_{d\in D}A_{d}: \varphi_{d,e}(x_{e})=x_{d}
\mbox{ for all }d,e\in D\mbox{ such that }
\, d\leq e\}.
\end{equation}
\noindent
In general, a pro-C$^{*}$-algebra is not a C$^{*}$-algebra
(see ``1.3. EXAMPLE" in \cite{Voiculescu}),
but the inverse limit of a special inverse system of C$^{*}$-algebras
is also a C$^{*}$-algebra as follows.
%
%
\begin{fact}
\label{fact:general}
Let $\{(A_{d},\varphi_{d,e}):d,e\in D\}$
be an inverse system of C$^{*}$-algebras.
We write $\varprojlim(A_{d},\varphi_{d,e})$ as
$\varprojlim A_{d}$ for simplicity of description.
\begin{enumerate}
\item
If $\varphi_{d,e}$ is injective for each $d,e$,
then
$\varprojlim A_{d}$ is a C$^{*}$-algebra.
\item
If $\{A_{n}:n\in {\bf N}\}$
is a sequence of inclusions of C$^{*}$-algebras
such that $A_{n}\supset A_{n+1}$ for each $n\in {\bf N}$,
then the inclusion map $\iota_{n}:A_{n+1}\hookrightarrow A_{n}$
induces an inverse system over $({\bf N},\leq)$ such that
$\varprojlim A_{n}$ is $*$-isomorphic onto $\bigcap_{n\geq 1}A_{n}$
as a C$^{*}$-algebra.
\item
If $A_{d}$ is unital and simple, and $f_{d,e}$ is unital for each $d,e\in D$,
then 
$\varprojlim A_{d}$ is a unital C$^{*}$-algebra.
\item
If there exists the maximal element $\omega$ of $D$,
then $\varprojlim A_{d}\cong A_{\omega}$.
\end{enumerate}
\end{fact}
%
%
\pr
Since an injective $*$-homomorphism is an isometry,
(i) holds from Definition \ref{defi:inverse}(ii)(c).
(ii) and (iii) hold from (i).
(iv) holds from (\ref{eqn:xain}).
\qedh

\noindent
When an inverse limit of C$^{*}$-algebras is a C$^{*}$-algebra,
it is not interesting as a pro-C$^{*}$-algebra, but
it does not mean the triviality of the inverse limit
as a C$^{*}$-algebra.

Next, we consider the inverse limit of Cuntz algebras in general 
setting.
(About inductive (or direct) limits of Cuntz algebras,
see \cite{GPW,LP,Rordam1993}.)
Since any Cuntz algebra is unital and simple,
the inverse limit of any inverse system of Cuntz algebras
by unital $*$-homomorphisms 
is a unital C$^{*}$-algebra from
Fact \ref{fact:general}(iii).
By Corollary \ref{cor:second} and Fact \ref{fact:general}(iv),
the following holds.
%
%
\begin{fact}
\label{fact:ica}
Let $(\hat{{\bf N}},\preceq)$ be as in Corollary \ref{cor:second}, and 
let $\{(\co{n(d)+1},\varphi_{d,e}):d,e\in D\}$
be an inverse system of Cuntz algebras
over a directed set $(D,\leq)$
such that $\{n(d)\in \hat{{\bf N}}:d\in D\}$.
We assume that $\varphi_{d,e}$ is unital for each $d,e$.
\begin{enumerate}
\item
The map
%
%
\begin{equation}
\label{eqn:fmap}
F:D\to \hat{{\bf N}};\quad d\mapsto F(d)\equiv n(d)
\end{equation}
is an ordered set homomorphism from
$(D,\leq)$ to $(\hat{{\bf N}},\preceq)$.
%
\item
If there exists the maximal element $\omega$ of $D$,
then $\varprojlim \co{n(d)+1}\cong \co{n(\omega)+1}$.
\end{enumerate}
\end{fact}

\noindent
Remark that $F$ in (\ref{eqn:fmap}) is not injective in general.
There are many endomorphisms of Cuntz algebras \cite{BJP, BJ1997,PE04}.

Let $(\hat{{\bf N}},\preceq)$ be as in Corollary 
\ref{cor:second}.
Define the order $\preceq_{c}$ on the
set $\{\co{n+1}:n\in \hat{{\bf N}}\}$
of all Cuntz algebras
as ``$A\preceq_{c}B$ if and only if ${\rm Hom}(A,B)\ne \emptyset$."
Then
$(\{\co{n+1}:n\in \hat{{\bf N}}\},\,\preceq_{c})$
is anti-isomorphic onto 
$(\hat{{\bf N}},\preceq)$ as an ordered set
with respect to the mapping $n\mapsto \co{n+1}$.

%
%
\ssft{An inverse system of Cuntz algebras}
\label{subsection:firstthree}
In this subsection,
we construct an example of inverse system of Cuntz algebras
as an answer to Question \ref{que:first}(i).
For the directed set $({\bf N},\preceq)$
in (\ref{eqn:precdef}),
define the inverse system $\{(R_{n},f_{n,m}):n,m\in {\bf N}\}$
of C$^{*}$-algebras as follows:
For $2\leq n<\infty$,
let $s_{1}^{(n)},\ldots,s_{n}^{(n)}$ denote 
the Cuntz generators of $\con$, that is, 
$(s_{i}^{(n)})^{*}s_{j}^{(n)}=\delta_{ij}I$ 
for $i,j=1,\ldots,n$
and 
$s_{1}^{(n)}(s_{1}^{(n)})^*+\cdots+s_{n}^{(n)}(s_{n}^{(n)})^*=I $.
For convenience,
rewrite $\co{n+1}$ as 
%
%
\begin{equation}
\label{eqn:pone}
R_{n}\equiv \co{n+1}\quad(n\in {\bf N}).
\end{equation}
This notation is reasonable with respect to
the $K$-theory of Cuntz algebras and 
Corollary \ref{cor:second}.
(About such a notation, see also \cite{ACE},
which is not related to the inclusions in Figure \ref{fig:three}.)
Remark that $R_{n}$
is generated by 
$s_{1}^{(n+1)},\ldots,s_{n+1}^{(n+1)}$ by definition.
When $n\preceq m$ and $n\ne m$,
define the $*$-homomorphism
$f_{n,m}$ from $R_{m}$ to $R_{n}$ by
%
%
\begin{equation}
\label{eqn:gencc}
\left\{
\begin{array}{rll}
f_{n,m}(s^{(m+1)}_{nl+i})\equiv
&(s_{n+1}^{(n+1)})^{l}s_{i}^{(n+1)}\quad&
\left(
\begin{array}{l}
l=0,1,\ldots,\frac{m}{n}-1,\\
i=1,\ldots,n
\end{array}
\right),\\
&\\
f_{n,m}(s^{(m+1)}_{m+1})\equiv &(s_{n+1}^{(n+1)})^{\frac{m}{n}}
\end{array}
\right.
\end{equation}
where $(s_{n+1}^{(n+1)})^{0}$ means the unit of $R_{n}$
for convenience,
and 
define
$f_{n,n}$ as the identity map $id_{R_{n}}$ on $R_{n}$.
A graphical explanation of (\ref{eqn:gencc})
will be given in $\S$ \ref{subsection:thirdtwo}.
Then the following holds.
%
%
\begin{Thm}
\label{Thm:inverseb}
Let $R_{n}$ and $f_{n,m}$ be as 
in (\ref{eqn:pone}) and (\ref{eqn:gencc}), respectively.
\begin{enumerate}
\item
The data $\{(R_{n},f_{n,m}):n,m\in {\bf N}\}$
is an inverse system of C$^{*}$-algebras
over the directed set $({\bf N},\preceq)$,
that is,
the relation
$f_{n,m}\circ f_{m,l}=f_{n,l}$ holds
when $n\preceq m\preceq l$.
\item
Let $\{s^{(\infty)}_{1},s^{(\infty)}_{2},\ldots\}$ denote
the Cuntz generators of $\coni$.
For $n\in {\bf N}$,
define the embedding $f_{n,\infty}$ of $\coni$ into $R_{n}$ by 
%
%
\begin{equation}
\label{eqn:eone}
f_{n,\infty}(s_{ln+i}^{(\infty)})\equiv (s_{n+1}^{(n+1)})^{l}s_{i}^{(n+1)}
\quad(i=1,\ldots,n,\,l\geq 0).
\end{equation}
Then $\{f_{n,\infty}:n\in {\bf N}\}$ satisfies
%
%
\begin{equation}
\label{eqn:unione}
f_{n,m}\circ f_{m,\infty}=f_{n,\infty}\quad\mbox{when }n\preceq m.
\end{equation}
\item
Every $f_{n,m}$ in (\ref{eqn:gencc}) and (\ref{eqn:eone})
is irreducible, 
where we state that $f\in {\rm Hom}(A,B)$ is irreducible if 
$f(A)^{'}\cap B={\bf C}I$.
\end{enumerate}
\end{Thm}

We illustrate relations of maps 
in Theorem \ref{Thm:inverseb}
as the following commutative diagrams
where we assume $n\preceq m$:

\noindent 
%
\def\setone{
\put(-200,0){$\coni$}
\put(200,200){$R_{m}=\co{m+1}$}
\put(200,-200){$R_{n}=\co{n+1}$}
\put(650,0){$R_{1}=\co{2}$}
\put(-100,50){\vector(2,1){250}}
\put(-100,-20){\vector(2,-1){250}}
\put(220,160){\vector(0,-1){280}}
\put(370,170){\vector(2,-1){250}}
\put(370,-150){\vector(2,1){250}}
\put(-50,140){$f_{m,\infty}$}
\put(-50,-120){$f_{n,\infty}$}
\put(240,0){$f_{n,m}$}
\put(500,140){$f_{1,m}$}
\put(500,-120){$f_{1,n}$}
\put(100,20){\rotatebox{90}{$\circlearrowright$}}
\put(400,20){$\circlearrowleft$}
}
\nolinebreak
%
%
%
\setlength{\unitlength}{.1mm}
\begin{picture}(1811,550)(-20,-230)
\thicklines
\put(-50,295){
\begin{minipage}[t]{2in}
\begin{fig}
\label{fig:first}
\end{fig}
\end{minipage}
}
\put(300,0){\setone}
\end{picture}

%
%
\begin{rem}
\label{rem:first}
{\rm
Every $f_{n,m}$ in (\ref{eqn:gencc}) and (\ref{eqn:eone})
 is unital and injective, but not surjective when $n\ne m$.
Such $*$-homomorphisms
were given in \cite{CK01,SE01}, hence they are not new,
but their relations of inverse system are new.
The essential part of 
Theorem \ref{Thm:inverseb}(i) 
is the construction of formulas in (\ref{eqn:gencc}).
The proof is given by simple algebraic calculation.
We make a point that 
Corollary \ref{cor:second} itself does not show
the existence of any inverse system of Cuntz algebras over the directed set
$({\bf N},\preceq)$.
Inversely, it is interesting question
whether the existence of such an inverse system of Cuntz algebras with unital
$*$-homomorphisms is shown only from 
Corollary \ref{cor:second} and general theory without use
of concrete construction of $*$-homomorphisms.
}
\end{rem}

For an example of $f_{n,m}$ in (\ref{eqn:gencc}),
the $*$-homomorphism $f_{1,2}:R_{2}(=\co{3})\to R_{1}(=\co{2})$ 
is given as follows:
%
%
\begin{equation}
\label{eqn:asada}
f_{1,2}(\hat{s}_{1})=s_{1},\quad 
f_{1,2}(\hat{s}_{2})=s_{2}s_{1},\quad f_{1,2}(\hat{s}_{3})=s_{2}s_{2}
\end{equation}
where $s_{1},s_{2}$ and $\hat{s}_{1},\hat{s}_{2},\hat{s}_{3}$
denote Cuntz generators of $\co{2}$ and $\co{3}$, respectively.
A similar $*$-homomorphism was given by Cuntz's original paper 
(\cite{Cuntz}, p183). 
For the first time,
the author knew (\ref{eqn:asada}) from Akira Asada who constructed $f_{1,2}$, 
and (\ref{eqn:gencc}) is a generalization of (\ref{eqn:asada}).

%
%
\ssft{Inverse limits of Cuntz algebras}
\label{subsection:firstfour}
In this subsection,
we show the inverse limit of the inverse system
$\{(R_{n},f_{n,m}):n,m\in {\bf N}\}$ in 
Theorem \ref{Thm:inverseb}(i),
which is the answer to Question \ref{que:first}(ii).
The problem is solved in slightly generalized 
setting.
%
%
\begin{Thm}
\label{Thm:inversec}
Let $\{(R_{n},f_{n,m}):n,m\in {\bf N}\}$ be as in 
Theorem \ref{Thm:inverseb}.
For  a directed subset $\Lambda$ of $({\bf N},\preceq)$,
let $\hatcola$ denote the inverse limit $\varprojlim_{\Lambda} R_{n}$ of 
the subsystem $\{(R_{n},f_{n,m}):n,m\in \Lambda\}$
of the inverse system $\{(R_{n},f_{n,m}):n,m\in {\bf N}\}$:
%
%
\begin{equation}
\label{eqn:hatcola}
\hatcola\equiv \varprojlim_{\Lambda} R_{n}.
\end{equation}
Then the following holds:
\begin{enumerate}
\item
Assume that $\Lambda$ is an infinite totally ordered subset of 
$({\bf N},\preceq)$
such that $\Lambda=\{n_{1},n_{2},\ldots\}$ and
$n_{1}\preceq n_{2}\preceq\cdots$.
For $\{f_{n,\infty}:n\in {\bf N}\}$ in (\ref{eqn:eone}),
define the $*$-homomorphism $\psi_{\Lambda}$ from $\coni$ to $\hatcola$
by
%
%
\begin{equation}
\label{eqn:psimap}
\psi_{\Lambda}(x)\equiv (f_{n_{1},\infty}(x),
f_{n_{2},\infty}(x),\ldots)\quad(x\in \coni)
\end{equation}
where $\hatcola$ is identified with the standard form
(\ref{eqn:xain}).
Then $\psi_{\Lambda}$ is a $*$-isomorphism such that
%
%
\begin{equation}
\label{eqn:unitwo}
\pi_{n}\circ \psi_{\Lambda}=f_{n,\infty}\quad(n\in \Lambda)
\end{equation}
where $\pi_{n}$ denotes the canonical homomorphism from $\hatcola$ to
$R_{n}$.
\item
For a directed subset $\Lambda$ of $({\bf N},\preceq)$,
the following holds:
%
%
\begin{equation}
\label{eqn:mainone}
\disp{\hatcola}\cong
\left\{
\begin{array}{ll}
\co{N+1}\quad &(\exists N=\max \Lambda),\\
\\
\coni \quad &(\mbox{otherwise}).
\end{array}
\right.
\end{equation}
Especially, when $\Lambda={\bf N}$, we obtain
%
%
\begin{equation}
\label{eqn:isom}
\varprojlim R_{n}\cong \coni.
\end{equation}
\end{enumerate}
\end{Thm}

\noindent
The proof of Theorem \ref{Thm:inversec} will be given 
in $\S$ \ref{section:second}.
The crucial part of the proof is the surjectivity of $\psi_{\Lambda}$ in
(\ref{eqn:psimap}).
From Theorem \ref{Thm:inversec}(i) for $\Lambda={\bf N}$,
the data
$(\coni,\{f_{n,\infty}\}_{n\in {\bf N}})$ is the inverse limit 
of $\{(R_{n},f_{n,m}):n,m\in {\bf N}\}$
in the sense of Definition \ref{defi:inverse}(ii).

As a counterview of Theorem \ref{Thm:inversec}(ii),
we can say that
$\coni$ can be written as an inverse limit
of $\con$'s.
From this and Definition \ref{defi:inverse}(ii)(b),
the following corollary immediately holds.
%
%
\begin{cor}
\label{cor:universal}
Let 
$\{f_{n,m}:n,m\in {\bf N}\}$ and
$\{f_{n,\infty}:n\in {\bf N}\}$ 
be as in (\ref{eqn:gencc}) and (\ref{eqn:eone}), respectively.
If a $*$-algebra $B$ and $*$-homomorphisms
$\{\eta_{n}\}_{n\in {\bf N}}$,
$\eta_{n}:B\to \co{n+1}$,
which satisfy $f_{n,m}\circ \eta_{m}=\eta_{n}$ when $n\preceq m$,
there exists a unique $*$-homomorphism $\psi$ from $B$ to $\coni$
such that $f_{n,\infty}\circ \psi=\eta_{n}$ for each $n\in {\bf N}$.
\end{cor}

\noindent
This is an essentially new universal property of $\coni$.
About other universal properties of $\coni$,
see Chapter 7 of \cite{RS}. 
%
%
\begin{rem}
\label{rem:firstb}
{\rm
\begin{enumerate}
\item
Here we discuss the choice of notations.
Both notations $\con$ and $R_{n}$
are useful in different situations.
The standard notations
$\co{2},\co{3},\ldots$ are used in the 
construction of 
C$^{*}$-bialgebra \cite{TS02}.
With respect to the standard notations of Cuntz algebras,
the order $\preceq$ in (\ref{eqn:precdef})
can be rewritten as follows:
A new order
$\ll$ on the set ${\bf N}_{\geq 2}\equiv \{2,3,4,\ldots\}$ 
is defined as
``$n\ll m$ if $n-1\preceq m-1$,"
or ``$\frac{m-1}{n-1}\in {\bf N}$."
Let $g_{n,m}\equiv f_{n-1,m-1}$ for $n,m\in {\bf N}_{\geq 2}$.
Then
$\{(\con,g_{n,m}): n,m\in {\bf N}_{\geq 2}\}$
is an inverse system over the directed set
$({\bf N}_{\geq 2},\ll)$
which is isomorphic to the inverse system 
$\{(R_{n},f_{n,m}):n,m\in {\bf N}\}$ over 
$({\bf N},\preceq)$.
By using $({\bf N}_{\geq 2},\ll)$,
the formula (\ref{eqn:firstlimit}) makes sense.
\item
In ``otherwise" in (\ref{eqn:mainone}),
we assume that $\Lambda$ is neither a cofinal 
nor a totally ordered
subset of 
$({\bf N},\preceq)$.
\item
We consider 
inverse systems of C$^{*}$-subalgebras of Cuntz algebras 
in Theorem \ref{Thm:inverseb}.
Let $\gamma^{(n+1)}$ and $\eta^{(n+1)}$
denote the $U(1)$-gauge action and 
the standard torus (=${\bf T}^{n+1}$)-action on $\co{n+1}$, respectively.
Let ${\cal A}_{n}$ and ${\cal C}_{n}$ denote the fixed-point 
subalgebras of $\co{n+1}$ with respect to 
$\gamma^{(n+1)}$ and $\eta^{(n+1)}$, respectively:
%
%
\begin{equation}
\label{eqn:fixedalgebra}
{\cal A}_{n}\equiv (\co{n+1})^{U(1)},\quad
{\cal C}_{n}\equiv (\co{n+1})^{{\bf T}^{n+1}}.
\end{equation}
Then we see that $f_{n,m}({\cal A}_{m})\not \subset {\cal A}_{n}$
even if $n\preceq m$, but
$f_{n,m}({\cal C}_{m})\subset {\cal C}_{n}$.
In this way,
$\{({\cal C}_{n},f_{n,m}|_{{\cal C}_{m}}):n,m\in {\bf N}\}$
is an inverse system of unital abelian C$^{*}$-algebras
over $({\bf N},\preceq)$.
Since $f_{n,m}|_{{\cal C}_m}$ is also injective and unital,
$\varprojlim
({\cal C}_{n},f_{n,m}|_{{\cal C}_{m}})$ is also a unital abelian C$^{*}$-algebra.
\item
As analogy with previous works 
of inductive limits \cite{Li,LP,Rordam1993},
the classification of inverse limits of Cuntz algebras is an interesting problem.
In the case of inductive limits,
classification theorems are given by using $K$-theory.
On the other hand,
it seems that
$K$-theory is no use for 
the case of inverse limits,
which will be shown in $\S$ \ref{subsection:firstfive}.
In order to consider the problem,
concrete examples are not sufficient yet. 
We will show other example of inverse system which limit
is not a Cuntz algebra in $\S$ \ref{subsection:thirdone}.
\item
In the proof of ``3.11 COROLLARY" in \cite{Cuntz1981},
$\coni$ is written as an inductive limit of 
C$^{*}$-subalgebras of $\con$'s as follows:
For $n\geq 1$,
let $C^{*}(s_{1}^{(n+1)},\ldots,s_{n}^{(n+1)})$
denote the C$^{*}$-subalgebra of $\co{n+1}$
generated by
$s_{1}^{(n+1)},\ldots,s_{n}^{(n+1)}$
without $s_{n+1}^{(n+1)}$.
Remark $C^{*}(s_{1}^{(n+1)},\ldots,s_{n}^{(n+1)})\not\cong \con$
because $K_{0}(C^{*}(s_{1}^{(n+1)},\ldots,s_{n}^{(n+1)}))\cong {\bf Z}$.
Then it is identified with
a $C^{*}$-subalgebra of 
$C^{*}(s_{1}^{(n+2)},\ldots,s_{n+1}^{(n+2)})$ 
by $\iota_{n}(s_{i}^{(n+1)})\equiv s_{i}^{(n+2)}$
for $i=1,\ldots,n$.
Then the inductive system 
$\{(C^{*}(s_{1}^{(n+1)},\ldots,s_{n}^{(n+1)}),\,\iota_{n}):n\geq 1\}$ 
over the directed set $({\bf N},\leq)$ is obtained:
%
%
\begin{equation}
\label{eqn:inclusionsix}
C^{*}(s_{1}^{(2)})\,
\subset \,
C^{*}(s_{1}^{(3)},s_{2}^{(3)})\,
\subset \,
C^{*}(s_{1}^{(4)},s_{2}^{(4)},s_{3}^{(4)})\,\subset\, \cdots,
\end{equation}
and its inductive limit is isomorphic onto $\coni$:
%
%
\begin{equation}
\label{eqn:limittwo}
\varinjlim C^{*}(s_{1}^{(n+1)},\ldots,s_{n}^{(n+1)})\cong \coni.
\end{equation}
This is quite a contrast to our result.
\end{enumerate}
}
\end{rem}

%
%
\ssft{Discontinuity of $K_{0}$}
\label{subsection:firstfive}
We discuss the (dis-) continuity of $K_{0}$-functor
of C$^{*}$-algebras
as an answer to Question \ref{que:first}(iii).
Related to the question of the continuity of $K_{0}$-functor
with respect to the inverse limit,
we could not find similar results in the standard textbooks
\cite{BlackadarK,JT,RLL,WO}.
In $\S$ 3 of \cite{Phillips1989},
the inverse limit for representable $K$-theory
is considered  for $\sigma$-C$^{*}$-algebras 
as the Milnor $\varprojlim^{1}$-sequence. 
However, it is assumed that all maps in the inverse system are surjective
(``3.2 THEOREM", \cite{Phillips1989}).
Hence it is no use for our example.
%
%
\sssft{Profinite groups}
\label{subsubsection:firstfiveone}
In order to discuss $K$-groups of inverse limits of Cuntz algebras,
we recall basic examples of profinite group 
(especially, they are procyclic groups (\cite{RZ}, $\S$ 2.7))
as follows.
Let $({\bf N},\preceq)$ be as in (\ref{eqn:precdef}).
For $n,m\in {\bf N}$,
if $n\preceq m$,
then the natural projection
from ${\bf Z}/m{\bf Z}$ onto ${\bf Z}/n{\bf Z}$
induces an inverse system
$\{{\bf Z}/n{\bf Z}:n\in {\bf N}\}$ of finite cyclic groups
(especially, they are rings)
over $({\bf N},\preceq)$.
It is well-known that
the inverse limit $\hat{{\bf Z}}$ of $\{{\bf Z}/n{\bf Z}:n\in {\bf N}\}$ 
is called the {\it Pr\"{u}fer ring}
\cite{Neukirch}:
%
%
\begin{equation}
\label{eqn:zhat}
\hat{{\bf Z}}\equiv \varprojlim {\bf Z}/n{\bf Z}.
\end{equation}
Remark $\hat{{\bf Z}}\not \cong {\bf Z}$.
For a fixed prime number $p$,
the subset $\{p^{n}:n\geq 1\}$ of ${\bf N}$ is a directed subset of 
$({\bf N},\preceq)$ which is not cofinal.
For the subsystem $\{{\bf Z}/p^{n}{\bf Z}:n\geq 1\}$
of $\{{\bf Z}/n{\bf Z}:n\in {\bf N}\}$,
the {\it pro-$p$ group} ${\bf Z}_{p}$
is defined as follows
\cite{RZ}:
%
%
\begin{equation}
\label{eqn:propgroup}
{\bf Z}_{p}\equiv \varprojlim {\bf Z}/p^{n}{\bf Z}.
\end{equation}
Remark that ${\bf Z}_{p}$ is identified with 
a uncountable proper subgroup of 
$\hat{{\bf Z}}$ (\cite{RZ}, Exercise 2.1.8).

%
%
\sssft{$K_{0}$-functor is not continuous with respect to
the inverse limit}
\label{subsubsection:firstfivetwo}
It is well-known that 
the $K_{0}$-functor of C$^{*}$-algebras
is continuous with respect to the inductive limit 
as the following sense \cite{BlackadarK}:
%
%
\begin{equation}
\label{eqn:continuous}
K_{0}(\varinjlim A_{n})\cong \varinjlim K_{0}(A_{n}).
\end{equation}
Since $\varinjlim A_{n}$ is always a C$^{*}$-algebra for any
inductive system,
(\ref{eqn:continuous}) holds for any
inductive system of C$^{*}$-algebra.
On the other hand,
the inverse limit of C$^{*}$-algebra is not always a C$^{*}$-algebra.
Hence (standard) 
$K$-groups can not be defined on a pro-C$^{*}$-algebra in general.
($K$-groups are generalized for $\sigma$-C$^{*}$-algebras
\cite{Phillips1989}, see also \cite{Cuntz2005,Weidner1,Weidner2}.)

From (\ref{eqn:firstlimit}),
the $K_{0}$-group of $\varprojlim\con$ is well-defined.
Then the following holds:
%
%
\begin{equation}
\label{eqn:kzero}
K_{0}(\varprojlim \con)\cong 
K_{0}(\coni)\cong {\bf Z}\not\cong
\hat{{\bf Z}}=\varprojlim {\bf Z}/n{\bf Z}\cong 
\varprojlim K_{0}(\con).
\end{equation}
%
This shows that 
$K_{0}$-functor is discontinuous
with respect to the inverse limit
even if the limit is a C$^{*}$-algebra.
Since $\hat{{\bf Z}}$ is the profinite completion
of ${\bf Z}$ \cite{Neukirch},
the profinite completion of $K_{0}(\varprojlim \con)$ 
coincides with $\varprojlim K_{0}(\con)$.
However 
the profinite completion of $K_{0}(\varprojlim \co{p^{n}+1})\cong K_{0}(\coni)$
does not coincide with $\varprojlim K_{0}(\co{p^{n}+1})\cong {\bf Z}_{p}$
from Theorem \ref{Thm:inversec}(ii) and (\ref{eqn:propgroup}).
Hence the profinite completion
does not always recover the continuity of the $K_0$-functor
with respect to the inverse limit.

In $\S$ \ref{section:second}, we will prove 
Theorem \ref{Thm:inverseb} and
Theorem \ref{Thm:inversec}.
In $\S$ \ref{section:third}, we will show examples.

%
%
\sftt{Proofs of main theorems}
\label{section:second}
In this section, we prove main theorems in $\S$ \ref{section:first}.
%
\ssft{Proof of Theorem \ref{Thm:inverseb}}
\label{subsection:secondone}
In this subsection,
we prove Theorem \ref{Thm:inverseb}.
By simple calculation, 
both (i) and (ii) are directly verified from (\ref{eqn:gencc}) and
(\ref{eqn:eone}).
In order to prove (iii),
we recall a lemma.
%
%
\begin{lem}
\label{lem:irreducible}
\begin{enumerate}
\item
Let $A$ and $B$ be unital C$^{*}$-algebras
and let $\rho$ be a unital $*$-homomorphism from $A$ to $B$.
If $B$ is simple and 
there exists an irreducible representation $\pi$ of $B$
such that
$\pi\circ \rho$ is also irreducible,
then $\rho$ is irreducible.
\item 
Let $\{s_{i}^{(n)}\}$ denote Cuntz generators 
of $\con$ for $2\leq n\leq \infty$.
Fix $i\in \{1,\ldots,n\}$.
Then there exists a unique state $\omega$ 
on $\con$ such that $\omega(s_{i}^{(n)})=1$.
Furthermore, $\omega$ is pure.
\end{enumerate}
\end{lem}
%
%
\pr
(i)
This holds from Proposition 3.1 in \cite{SE01}.

\noindent
(ii)
The uniqueness holds by Cuntz relations.
Let $({\cal H}_{\omega},\pi_{\omega},\Omega_{\omega})$
denote the GNS triple by $\omega$.
By definition, we see 
$\pi_{\omega}(s_{i}^{(n)})\Omega_{\omega}=\Omega_{\omega}$.
From \cite{BJ,DaPi2,DaPi3},
such a cyclic representation $\pi_{\omega}$ exists and is irreducible.
Hence $\omega$ is pure.
\qedh

Let $\omega_{n}$ denote the state on $R_n=\co{n+1}$
such that 
%
%
\begin{equation}
\label{eqn:stateone}
\omega_{n}(s_{1}^{(n+1)})=1.
\end{equation}
From Lemma \ref{lem:irreducible}(ii),
the GNS representation $\pi_{n}$ by $\omega_{n}$
is irreducible.
When $n\preceq m$,
we can verify 
$\omega_{n}\circ f_{n,m}=\omega_{m}$
because $f_{n,m}(s_{1}^{(m+1)})=s_{1}^{(n+1)}$.
Hence $\pi_{n}\circ f_{n,m}$ is unitarily equivalent to $\pi_{m}$,
and $\pi_{n}\circ f_{n,m}$ is also irreducible.
From this and Lemma \ref{lem:irreducible}(i),
Theorem \ref{Thm:inverseb}(iii) is proved.


%
%
\ssft{Proof of Theorem \ref{Thm:inversec}}
\label{subsection:secondtwo}
In this subsection,
we prove Theorem \ref{Thm:inversec}
 except
a certain equality of C$^{*}$-subalgebras,
which will be proved in $\S$ \ref{subsection:secondthree}.

In order to reduce the problem,
we show a lemma.
A subset $E$ of a directed set $(D,\leq)$
is {\it cofinal} if
$\{e\in D:e\geq d\}\cap E\ne\emptyset$ for any $d\in D$.
In this case,
it is known that $\varprojlim_{E}A_{d}\cong \varprojlim_{D}A_{d}$
for the subsystem 
$\{(A_{d},\varphi_{d,e}):d,e\in E\}$
of the inverse system 
$\{(A_{d},\varphi_{d,e}):d,e\in D\}$ over $(D,\leq)$.
Then the following lemma holds.
%
%
\begin{lem}
\label{lem:engine}
For any countable directed set $(D,\leq)$,
there exists a totally ordered cofinal subset $D_{0}$ of $D$.
\end{lem}
%
%
\pr
If the maximal element $\omega$ of $D$ exists,
then let $D_0 \equiv \{\omega\}$.
If not,
let $D=\{x_{1},x_{2},\ldots\}$,
where we do not assume $x_{1}\leq x_{2}\leq \cdots$.
Then we can inductively construct
a subsequence $y_{1},y_{2},\ldots$ of $x_{1},x_{2},\ldots$
as follows:
Let $y_{1}=x_{1}$.
For $n\geq 2$,
there always exists $z\in D$ such that
$y_{n-1}\leq z$ and $x_{n}\leq z$.
Choose such an element $z$ and 
define $y_{n}\equiv z$.
Then $D_0\equiv \{y_{1},y_{2},\ldots\}$ is a totally ordered cofinal subset of $D$.
\qedh 

\noindent
For example,
$\{n!:n\in {\bf N}\}$ is a totally ordered cofinal 
subset of $({\bf N},\preceq)$.

\ww
{\it Proof of Theorem \ref{Thm:inversec}.}
(i)
Here we will prove the statement except a certain equality.
From Fact \ref{fact:general}(iii), $\hatcola$ is a C$^{*}$-algebra.
By the standard construction 
in (\ref{eqn:xain}),
$\hatcola$ is given as follows:
%
%
\begin{equation}
\label{eqn:inverse}
\{(x_{n_1},x_{n_2},\ldots)\in  
\prod_{k\geq 1}R_{n_k}:f_{n_k,n_l}(x_{n_l})=x_{n_k}
\mbox{ for each }k\leq l\}.
\end{equation}

Define the set $\{Q_{n}:1\leq n\leq \infty\}$ 
of C$^{*}$-subalgebras of $R_{1}$ by
%
%
\begin{equation}
\label{eqn:qalgebra}
Q_{n}\equiv f_{1,n}(R_{n})\quad(1\leq n<\infty),\quad
Q_{\infty}\equiv f_{1,\infty}(\coni).
\end{equation}
Then we see that
$Q_{1}=R_{1}=\co{2}$,
$Q_{n}\cong R_{n}=\co{n+1}$ when $1\leq n<\infty$,
and 
$Q_{m}\subset Q_{n}$ when $n\preceq m$ from Theorem \ref{Thm:inverseb}(i).
On the other hand,
from (\ref{eqn:unione}),
$Q_{\infty}\subset Q_{n}$ for each $n\geq 1$.
Since $n_{1}\preceq n_{2}\preceq\cdots$,
we obtain the following unital inclusions
%
%
\begin{equation}
\label{eqn:inclusionsthree}
\co{2}=Q_{1}\supset Q_{n_{1}}\supset Q_{n_{2}}\supset Q_{n_{3}}\supset \cdots
\supset Q_{\infty}.
\end{equation}

Define 
$\hat{\pi}_{1}\equiv f_{1,n_{1}}\circ \pi_{n_{1}}$.
Then we see that 
%
%
\begin{equation}
\label{eqn:equalityone}
\hat{\pi}_{1}(\hatcola)=\bigcap_{n\in\Lambda} Q_{n}.
\end{equation}
Since $\hat{\pi}_{1}$ is injective,
$\hatcola$ and $\bigcap_{n\in\Lambda} Q_{n}$
are $*$-isomorphic, and the map
$f_{1,\infty}:\coni\to \bigcap_{n\in\Lambda}Q_{n}$
is well-defined.
In consequence, we see that
the following diagram is commutative:

\def\figures{
\put(0,0){$\coni$}
\put(100,15){\vector(1,0){100}}
\put(250,0){$\hatcola$}
\put(270,-30){\vector(0,-1){100}}
\put(240,-200){$\disp{\bigcap_{n\in\Lambda} Q_{n}}$}
\put(60,-40){\vector(1,-1){130}}
\put(140,30){$\psi_{\Lambda}$}
\put(150,-40){\rotatebox{-90}{{\bf $\circlearrowright$}}}
\put(130,-100){$f_{1,\infty}$}
\put(300,-90){$\hat{\pi}_{1}$}
}

%
%
%
\begin{fig}
\label{fig:onelast}
\quad\\
\thicklines
\setlength{\unitlength}{.15mm}
\begin{picture}(1000,250)(-200,-220)
\put(0,0){\figures}
\end{picture}
\end{fig}

In order to prove the bijectivity of $\psi_{\Lambda}$,
it is sufficient to show the bijectivity of $f_{1,\infty}$.
Since $f_{1,\infty}$ is an injective $*$-homomorphism,
it is sufficient to show that 
%
%
\begin{equation}
\label{eqn:identity}
Q_{\infty}=\bigcap_{n\in\Lambda} Q_{n}. 
\end{equation}
We will prove (\ref{eqn:identity})
in $\S$ \ref{subsection:secondthree}.

\noindent
(ii)
When there exists the maximal element of $\Lambda$,
the statement holds from Fact \ref{fact:ica}(ii).
Assume that $\Lambda$ has no maximal element.
In this case,
it is sufficient to assume the condition in (i)
for $\Lambda$ by Lemma \ref{lem:engine}.
Hence the statement holds from (i).
\qedh

%
%
\ssft{Proof of (\ref{eqn:identity})}
\label{subsection:secondthree}
From the proof of Theorem \ref{Thm:inversec}(i),
the problem is reduced to the relation (\ref{eqn:identity}) 
among C$^{*}$-subalgebras 
$\{Q_{n}:1\leq n\leq \infty\}$
of $Q_{1}(=R_{1}=\co{2})$ in (\ref{eqn:qalgebra}).
In this subsection, we prove (\ref{eqn:identity}).
%
%
\sssft{Inclusions of free subsemigroups of $Q_1$}
\label{subsubsection:secondthreeone}
In this subsubsection,
we consider free subsemigroups of $Q_1$ and their relations.
We rewrite the Cuntz generators of $Q_{1}$ as  $t_{1},t_{2}$ here.
Let $S(X)$ denote the subsemigroup of $Q_1$
generated by a subset $X$ of $Q_{1}$.
Define subsemigroups ${\cal K}_{n},{\cal L}_{n}$ of $Q_1$ as 
%
%
\begin{equation}
\label{eqn:subsemigroups}
\left\{
\begin{array}{rl}
{\cal K}_{n}\equiv &S(\{t_{2}^{n}\})\quad(1\leq n<\infty),\\
\\
{\cal L}_{1}\equiv &S(\{t_{1},t_{2}\}),\\
\\
{\cal L}_{n}\equiv &S(\{t_{1},t_{2}t_{1},\ldots,t_{2}^{n-1}t_{1},t_{2}^{n}\})
\quad(2\leq n<\infty),\\
\\
{\cal L}_{\infty}\equiv &S(\{t_{1},t_{2}^{m}t_{1}:m\geq 1\}).
\end{array}
\right.
\end{equation}
For each $n\geq 1$, ${\cal K}_{n}$ is abelian,
 and ${\cal L}_{n}$ 
is a (non-unital) free semigroup of rank $n+1$
\cite{Howie}.
Both ${\cal K}_{n}$ and ${\cal L}_{n}$ are subsemigroups of $Q_{n}$.
If $m\preceq n$, then ${\cal K}_{n}\subset {\cal K}_{m}$ and
${\cal L}_{n}\subset {\cal L}_{m}$.
For each $n\geq 1$, ${\cal L}_{n}\supset {\cal L}_{\infty}$. 
Furthermore, we see that
%
%
\begin{equation}
\label{eqn:concrete}
{\cal K}_{n}=\{t_{2}^{n},t_{2}^{2n},t_{2}^{3n},\ldots\},\quad
{\cal L}_{\infty}=\{xt_{1}:x\in {\cal L}_{1}\}={\cal L}_{1}t_1.
\end{equation}
%
Hence ${\cal K}_{n}\cap {\cal L}_{\infty}=\emptyset$.

Since ${\cal L}_{1}={\cal L}_{1}t_{1}\,\sqcup\, {\cal L}_{1}t_{2}$
with respect to the ending of each word,
the decomposition
${\cal L}_{n}
=
({\cal L}_{n}\cap{\cal L}_1 t_{1})
\sqcup 
({\cal L}_{n}\cap {\cal L}_{1}t_{2})$ holds.
Let $Y_{n}\equiv \{u,xu:x\in {\cal L}_{\infty},\,u\in {\cal K}_{n}\}$.
Then ${\cal L}_{n}\cap {\cal L}_1 t_{2}=Y_{n}$
and
${\cal L}_{n}\cap {\cal L}_1 t_{1}={\cal L}_{\infty}$.
Hence the following decomposition into disjoint subsets holds:
%
%
\begin{equation}
\label{eqn:decolast}
{\cal L}_{n}={\cal L}_{\infty}\sqcup Y_{n}\quad(2\leq n<\infty).
\end{equation}
%

%
%
\sssft{Decomposition of algebras into linear subspaces}
\label{subsubsection:secondthreetwo}
Let ${\cal K}_{n},{\cal L}_{n}$ be as in (\ref{eqn:subsemigroups}).
From definitions of $Q_{n}$
and $\{f_{n,m},f_{n,\infty}:n,m\in {\bf N}\}$,
we see that
%
%
\begin{equation}
\label{eqn:relativeb}
Q_{n}=C^{*}\langle {\cal L}_n \rangle
\quad(1\leq n\leq \infty)
\end{equation}
where $C^*\langle X\rangle$ denote the C$^{*}$-subalgebra of $Q_1$
generated by a subset $X$ of $Q_1$.
As closed linear subspaces of $Q_1$,
$Q_{n}$'s can be written as follows:
%
%
\begin{equation}
\label{eqn:relativec}
\left\{
\begin{array}{l}
Q_{n}=\overline{
{\rm Lin}\langle\{xy^{*}:x,y\in {\cal L}_{n}\}\rangle}
\quad(1\leq n<\infty),\\
\\
Q_{\infty}=
\overline{
{\rm Lin}\langle\{
I,xy^{*} :x,y\in {\cal L}_{\infty}\}\rangle}.
\end{array}
\right.
\end{equation}
Remark that 
vectors in each generating set in (\ref{eqn:relativec})
are not always linearly independent because of Cuntz relations.
%
%
%
\begin{lem}
\label{lem:decompositionlast}
\begin{enumerate}
\item
For $n\geq 1$,
$(t_{2})^{n}(t_{2}^{*})^{n}
=I-\sum_{k=0}^{n-1}(t_{2})^{k}t_{1}t_{1}^{*}(t_{2}^{*})^{k}$.
\item
Define closed linear subspaces $V_{n}$ and $V_{n}^{*}$ of $Q_{n}$ by
%
%
\begin{equation}
\label{eqn:vnone}
V_{n}\equiv\overline{{\rm Lin}\langle\{u,xu,xuy^{*}:
x,y\in {\cal L}_{\infty},\,u\in {\cal K}_{n}\}\rangle},\quad 
 V_{n}^{*}\equiv \{x^{*}:x\in  V_{n}\}.
\end{equation}
Then 
$V_{n}\subset V_{m}$ and
$V_{n}^{*}\subset V_{m}^{*}$ 
when $m\preceq n$,
$ V_{n}\cap  V_{n}^{*}=V_{n}\cap Q_{\infty}=V_{n}^{*}\cap Q_{\infty}=\{0\}$.
\item
For $u,v\in {\cal K}_{n}$ and $x,y\in {\cal L}_{\infty}$,
$xuv^{*}y^{*}\in
Q_{\infty}\oplus  V_{n}\oplus  V_{n}^{*}$.
\item
For each $n\geq 1$,
the following decomposition of $Q_{n}$ into closed linear subspaces holds:
%
%
\begin{equation}
\label{eqn:deconew}
Q_{n}=Q_{\infty}\oplus  V_{n}\oplus  V_{n}^{*}.
\end{equation}
\item
For $\Lambda$ in Theorem \ref{Thm:inversec}(i),
$\bigcap_{n\in\Lambda} V_{n}
=\bigcap_{n\in\Lambda} V_{n}^{*}=\{0\}$.
\end{enumerate}
\end{lem}
%
%
\pr
(i)
By the Cuntz relations of $Q_{1}=\co{2}$,
the statement holds.

\noindent
(ii)
By definition, the statement holds.

\noindent
(iii)
From (ii),
$Q_{\infty}\oplus  V_{n}\oplus  V_{n}^{*}$
makes sense as a subspace of $Q_n$.
If $u=v$, then
$xuv^{*}y^{*}\in Q_{\infty}$ from (i) and (\ref{eqn:relativec}).
Furthermore,
from (i),
$t_{2}^{n(l+k)}(t_{2}^{nk})^{*}=t_{2}^{nl}-t_{2}^{nl}(\sum_{j=0}^{nk-1}
t_{2}^{j}t_{1}t_{1}^{*}(t_{2}^*)^{j})\in  V_{n}\oplus Q_{\infty}$
for $l\geq 1$.
Hence the statement also holds for the case $u\ne v$.

\noindent
(iv)
From (ii),
$Q_{n}\supset Q_{\infty}\oplus  V_{n}\oplus  V_{n}^{*}$.
Since $Q_{n}$ is the closure of the linear space
spanned by the set
%
%
\begin{equation}
\label{eqn:mn}
{\cal M}_{n}\equiv \{x,x^{*},xy^{*}:x,y\in {\cal L}_{n}\}
\quad(1\leq n<\infty),
\end{equation}
it is sufficient to show that
${\cal M}_{n}$ is a subset of $Q_{\infty}\oplus  V_{n}\oplus  V_{n}^{*}$.
From (\ref{eqn:decolast}),
${\cal M}_{n}$ is decomposed into the disjoint union as follows:
%
%
\begin{equation}
\label{eqn:calm}
\begin{array}{rl}
{\cal M}_{n}
= &\{x,x^{*},xy^{*}:x,y\in {\cal L}_{\infty}\sqcup Y_{n}\}
={\cal M}_{n,1}\sqcup {\cal M}_{n,2},\\
\\
{\cal M}_{n,1}\equiv  &
\{x,x^{*},xy^{*}:x,y\in {\cal L}_{\infty}\},\\
\\
{\cal M}_{n,2}\equiv &
\{u,u^{*},uv^{*},xu^{*},ux^{*}:x\in {\cal L}_{\infty},\,u,v\in Y_{n}\}.\\
\end{array}
\end{equation}
Since
$Q_{\infty}=\overline{{\rm Lin}\langle {\cal M}_{n,1}\cup\{I\}\rangle}$,
it is sufficient to show that
${\cal M}_{n,2}$ is a subset of 
$Q_{\infty}\oplus  V_{n}\oplus  V_{n}^{*}$.
By the definition of $Y_{n}$ in (\ref{eqn:decolast}),
%
%
\begin{equation}
\label{eqn:wsettwo}
\begin{array}{rl}
{\cal M}_{n,2}=& {\cal M}_{n,2,1}\sqcup {\cal M}_{n,2,2},\\
\\
{\cal M}_{n,2,1}\equiv  &
\{u,u^{*},xu,(xu)^{*}:x\in {\cal L}_{\infty},\,u\in {\cal K}_{n}\},\\
\\
{\cal M}_{n,2,2}\equiv &\{xuy^{*}, xu^{*}y^{*},
xuv^{*}y^{*}:x,y\in {\cal L}_{\infty},u,v\in {\cal K}_{n}\}.\\
\end{array}
\end{equation}
Then 
${\cal M}_{n,2,1}\subset V_{n}\oplus  V_{n}^{*}$ by definition.
From (iii),
${\cal M}_{n,2,2}\subset Q_{\infty}\oplus  V_{n}\oplus  V_{n}^{*}$.
Hence 
${\cal M}_{n,2}\subset Q_{\infty}\oplus  V_{n}\oplus  V_{n}^{*}$.

\noindent
(v)
By definition,
%
%
\begin{equation}
\label{eqn:vdecompose}
V_{n}=
\overline{{\rm Lin}\langle {\cal K}_{n}\rangle }
\oplus 
\overline{{\rm Lin}\langle {\cal L}_{\infty}\cdot {\cal K}_{n}\rangle }
\oplus 
\overline{{\rm Lin}\langle {\cal L}_{\infty}\cdot {\cal K}_{n}\cdot
{\cal L}_{\infty}^{*}\rangle}
\end{equation}
where
${\cal L}_{\infty}^{*}\equiv \{x^{*}:x\in {\cal L}_{\infty}\}$.
By assumption, $\Lambda$ is an infinite subset of ${\bf N}$.
Hence $\bigcap_{n\in \Lambda}{\cal K}_{n}=\emptyset$.
From this,
(\ref{eqn:vdecompose}) and (i),
%
%
\begin{equation}
\label{eqn:vdecomposetwo}
\bigcap_{n\in\Lambda} V_{n}=
\bigcap_{n\in\Lambda}\overline{{\rm Lin}\langle {\cal K}_{n}\rangle }
\oplus 
\bigcap_{n\in\Lambda}\overline{{\rm Lin}\langle {\cal L}_{\infty}\cdot {\cal K}_{n}\rangle }
\oplus 
\bigcap_{n\in\Lambda}\overline{{\rm Lin}\langle {\cal L}_{\infty}\cdot {\cal K}_{n}\cdot
{\cal L}_{\infty}^{*}\rangle}=\{0\}.
\end{equation}
In a similar way, we obtain
$\bigcap_{n\in\Lambda} V_{n}^{*}=\{0\}$.
Hence the statement holds.
\qedh

From Lemma \ref{lem:decompositionlast}(iv) and (v),
%
%
\begin{equation}
\label{eqn:intersect}
\bigcap _{n\in\Lambda}Q_{n}
=Q_{\infty}\oplus 
(\bigcap_{n\in \Lambda} V_{n})
\oplus 
(\bigcap_{n\in \Lambda} V_{n}^{*})
=Q_{\infty}.
\end{equation}
Hence (\ref{eqn:identity}) is proved.

%
%
\sftt{Examples}
\label{section:third}
In order to explain 
theorems in $\S$ \ref{section:first} more,
we show examples in this section.

%
%
\ssft{Other inverse system}
\label{subsection:thirdone}
In this subsection, we show other example of inverse system
of Cuntz algebras.
Let $s_{1}^{(n)},\ldots,s_{n}^{(n)}$ denote
the Cuntz generators of $\con$.
Fix an integer $r\geq 2$.
Let $r_{n}\equiv r^{2^{n-1}}$ and rewrite $\co{r_{n}}$ as
%
%
\begin{equation}
\label{eqn:arndefi}
A_{r,n}\equiv \co{r_n}
\quad(n\geq 1).
\end{equation}
Then $r_{n}^{2}=r_{n+1}$ for $n\geq 1$.
Define the $*$-homomorphism
$q_{n}$ from $A_{r,n+1}$ to $A_{r,n}$ by
%
%
\begin{equation}
\label{eqn:qmap}
q_{n}(s_{r_n(i-1)+j}^{(r_{n+1})})\equiv 
s_{i}^{(r_n)}s_{j}^{(r_n)}
\quad(i,j=1,\ldots,r_n).
\end{equation}
Then
$\{(A_{r,n},q_{n}):n\geq 1\}$
is an inverse system of Cuntz algebras
over the directed set $({\bf N},\leq)$.

Remark $A_{r,n}=\co{(r_{n}-1)+1}$ for $n\geq 1$.
Here we verify Fact \ref{fact:ica}(i) 
for the sequence $\{r_{n}-1:n\geq 1\}$.
Define the map
$F$ from $({\bf N},\leq)$ to $({\bf N},\preceq)$ by
%
%
\begin{equation}
\label{eqn:exampletwo}
F(n)\equiv r_n-1\quad(n\in {\bf N}).
\end{equation}
Then 
$F(n)=r_n-1\preceq (r^{2^{n-1}}-1)(r^{2^{n-1}}+1)=r^{2^{n}}-1=F(n+1)$.
Exactly, the map $F$ satisfies the statement in Fact \ref{fact:ica}(i).

%
%
\begin{prop}
\label{prop:var}
The inverse limit 
$\varprojlim_{n} A_{r,n}$ 
of $\{(A_{r,n},q_{n}):n\in {\bf N}\}$
is $*$-isomorphic onto
the uniformly hyperfinite algebra $UHF_{r}$ 
of the Glimm type $\{r^{l}:l\geq 1\}$ \cite{Glimm}:
%
%
\begin{equation}
\label{eqn:glimm}
\varprojlim_{n} A_{r,n}\cong UHF_{r}.
\end{equation}
\end{prop}
%
%
\pr
Let $\gamma$ denote the $U(1)$-gauge action on
$\co{r}=\co{r_{1}}=A_{r,1}$.
For $l\in {\bf Z}$,
define
%
%
\begin{equation}
\label{eqn:arn}
A_{r,1}^{(l)}
\equiv \{x\in A_{r,1}:\gamma_{z}(x)=z^{l}x\mbox{ for all }z\in U(1)\}.
\end{equation}
Rewrite the Cuntz generators of $\co{r}=\co{r_{1}}$ as
$t_{1},\ldots,t_{r}$.
By identifying $A_{r,n}$
with the C$^{*}$-subalgebra
$(q_{1}\circ q_{2}\circ \cdots\circ q_{n-1})(A_{r,n})$ of $A_{r,1}$,
$A_{r,n}$'s are rewritten as follows:
%
%
\begin{equation}
\label{eqn:arnb}
A_{r,n}=C^{*}\langle\{t_{J}:J\in\{1,\ldots,r\}^{2^{n-1}}\}\rangle\quad(n\geq 1)
\end{equation}
where $t_{J}\equiv t_{j_{1}}\cdots t_{j_{m}}$
for $J=(j_{1},\ldots,j_{m})$.
Then we obtain the following unital (rapidly decreasing) inclusions: 
%
%
\begin{equation}
\label{eqn:ainclusion}
A_{r,1}\supset A_{r,2}\supset A_{r,3}\supset\cdots.
\end{equation}
Define
%
%
\begin{equation}
\label{eqn:subtraction}
A_{r,n}^{(l)}
=A_{r,n}\cap A_{r,1}^{(l)}\quad(l\in {\bf Z},\,n\geq 1).
\end{equation}
Then 
$A_{r,n+1}^{(l)}\subset A_{r,n}^{(l)}$ for each $n,l$,
and
$A_{r,n}^{(l)}$
is the closure of 
${\rm Lin}\langle
\{t_{J}t_{K}^{*}:J,K\in \bigcup_{a\geq 1}\{1,\ldots,r\}^{a\times 2^{n-1}},\,
|J|-|K|=l\}\rangle$.
When $l\ne 0$,
$\bigcap_{n\geq 1}A_{r,n}^{(l)}=\{0\}$ 
because $A_{r,n}^{(l)}=\{0\}$ if $2^{n-1}\nmid  l$.
Since $A_{r,n}=\bigoplus_{l\in {\bf Z}}A_{r,n}^{(l)}$
and $A_{r,n+1}^{(0)}= A_{r,n}^{(0)}$ for each $n$,
$\varprojlim_n A_{r,n}\cong \bigcap_{n\geq 1}A_{r,n}=A_{r,1}^{(0)}\cong
UHF_{r}$.
\qedh

Since $K_{0}(UHF_{r})$ is the group ${\bf Z}_{(r^{\infty})}\subset {\bf Q}$
of all rational numbers whose denominators divide 
the generalized integer $r^{\infty}$ 
($\S$ 7.5, \cite{BlackadarK}),
we see that
%
%
\begin{equation}
\label{eqn:violation}
\begin{array}{rl}
K_{0}(\varprojlim A_{r,n})\cong &
K_{0}(UHF_{r}) \cong {\bf Z}_{(r^{\infty})},\\
\\
\varprojlim K_{0}(A_{r,n})
=&\varprojlim K_{0}(\co{r_{n}})
\cong \varprojlim {\bf Z}/(r_{n}-1){\bf Z}.\\
\end{array}
\end{equation}
%
The former is countable, but the latter is not.
This case also shows that 
the $K_{0}$-functor is  discontinuous
with respect to the inverse limit.

%
%
\ssft{Illustrations of embeddings}
\label{subsection:thirdtwo}
In this subsection,
we illustrate embeddings in (\ref{eqn:gencc})
by using decompositions of Hilbert spaces.
When $\con$ acts on a Hilbert space,
by identifying a generator $s_{i}$
with the range of $s_{i}$,
the following illustration is helpful in understanding $s_{1},\ldots,s_{n}$:

%
%
\def\firstbox{
\put(0,0){\line(1,0){500}}
\put(0,30){\line(1,0){500}}
\put(0,0){\line(0,1){30}}
\put(500,0){\line(0,1){30}}
\put(240,10){${\huge {\cal H}}$}
}
\def\secondbox{
\put(0,0){\line(1,0){500}}
\put(0,30){\line(1,0){500}}
\put(0,0){\line(0,1){30}}
\put(100,0){\line(0,1){30}}
\put(200,0){\line(0,1){30}}
\put(300,0){\line(0,1){30}}
\put(140,10){$\cdots$}
\put(340,10){$\cdots$}
\put(400,0){\line(0,1){30}}
\put(500,0){\line(0,1){30}}
\put(40,10){${\huge s_{1}{\cal H}}$}
\put(240,10){${\huge s_{i}{\cal H}}$}
\put(440,10){${\huge s_{n}{\cal H}}$}
}
\def\cross{
\thinlines
\put(225,60){$s_{i}$}
\put(245,60){$\downarrow$}
\qbezier[200](500,100)(400,65)(300,30)
\qbezier[200](0,100)(100,65)(200,30)
}
%
%
%

\noindent
\setlength{\unitlength}{.22mm}
\begin{picture}(550,220)(-25,-30)
\put(-30,170){\begin{minipage}[t]{2in}
%
\begin{fig}
\label{fig:bone}
\end{fig}
\end{minipage}
}
\thicklines
\put(0,100){\firstbox}
\put(0,0){\secondbox}
\put(0,0){\cross}
\end{picture}

Recall that
$R_{1}=\co{2}$, $R_{2}=\co{3}$,
$R_{4}=\co{5}$. Then
$f_{1,2},f_{1,4},f_{2,4}$ in (\ref{eqn:gencc})
are given as follows:
\\
$f_{1,2}:R_{2}\to R_{1}$;
%
%
\begin{equation}
\label{eqn:fff}
f_{1,2}(s_{1}^{(3)})=s_{1}^{(2)},\quad 
f_{1,2}(s_{2}^{(3)})=s_{2}^{(2)}s_{1}^{(2)},\quad 
f_{1,2}(s_{3}^{(3)})=(s_{2}^{(2)})^{2}.
\end{equation}
$f_{1,4}:R_{4}\to R_{1}$;
%
%
\begin{equation}
\label{eqn:ffftwo}
f_{1,4}(s_{1}^{(5)})=s_{1}^{(2)},\quad
f_{1,4}(s_{2}^{(5)})=s_{2}^{(2)}s_{1}^{(2)},\quad
f_{1,4}(s_{3}^{(5)})=(s_{2}^{(2)})^{2}s_{1}^{(2)},
\end{equation}
%
%
%
\begin{equation}
\label{eqn:fffthree}
f_{1,4}(s_{4}^{(5)})=(s_{2}^{(2)})^{3}s_{1}^{(2)},\quad
f_{1,4}(s_{4}^{(5)})=(s_{2}^{(2)})^{4}.
\end{equation}
$f_{2,4}:R_{4}\to R_{2}$;
%
%
\begin{equation}
\label{eqn:ffffour}
f_{2,4}(s_{1}^{(5)})=s_{1}^{(3)},\quad
f_{2,4}(s_{2}^{(5)})=s_{2}^{(3)},\quad
f_{2,4}(s_{3}^{(5)})=s_{3}^{(3)}s_{1}^{(3)},
\end{equation}
%
%
%
\begin{equation}
\label{eqn:ffffive}
f_{2,4}(s_{4}^{(5)})=s_{3}^{(3)}s_{2}^{(3)},\quad
f_{2,4}(s_{5}^{(5)})=(s_{3}^{(3)})^{2}.
\end{equation}
From these,
we can directly verify the identity
$f_{1,2}\circ f_{2,4}=f_{1,4}$.

From Figure \ref{fig:bone},
inclusions $\co{2}\supset \co{3}\supset \co{5}$ 
by embeddings $f_{1,2}$ and $f_{2,4}$
are illustrated as follows:

\def\kbox#1#2{
\put(0,0){\line(1,0){#1}}
\put(0,#2){\line(1,0){#1}}
\put(0,0){\line(0,1){#2}}
\put(#1,0){\line(0,1){#2}}
}
\def\otwo{
\put(0,0){\kbox{500}{30}}
\put(250,0){\line(0,1){30}}
\multiput(375,0)(0,5){6}{\line(0,1){3}}
}
\def\othree{
\put(0,0){\kbox{500}{30}}
\put(170,0){\line(0,1){30}}
\put(340,0){\line(0,1){30}}
\multiput(393,0)(0,5){6}{\line(0,1){3}}
\multiput(450,0)(0,5){6}{\line(0,1){3}}
}
\def\ofive{
\put(0,0){\kbox{500}{30}}
\put(100,0){\line(0,1){30}}
\put(200,0){\line(0,1){30}}
\put(300,0){\line(0,1){30}}
\put(400,0){\line(0,1){30}}
}
\def\mainfig{
\put(0,100){\otwo}
\put(0,0){\othree}
\put(0,-100){\ofive}
\thinlines
\path(0,-70)(0,100)
\path(500,-70)(500,100)
\path(100,-70)(170,0)
\path(200,-70)(340,0)
\path(300,-70)(393,0)
\path(400,-70)(450,0)
\path(170,30)(250,100)
\path(340,30)(375,100)
\put(-30,10){$\co{3}$}
\put(-30,110){$\co{2}$}
\put(-30,-90){$\co{5}$}
\put(70,40){$s^{(3)}_{1}$}
\put(240,40){$s^{(3)}_{2}$}
\put(400,40){$s^{(3)}_{3}$}
\put(120,140){$s^{(2)}_{1}$}
\put(370,140){$s^{(2)}_{2}$}
\put(300,107){$s^{(2)}_{2}s^{(2)}_{1}$}
\put(420,107){$s^{(2)}_{2}s^{(2)}_{2}$}
\put(50,-60){$s^{(5)}_{1}$}
\put(130,-60){$s^{(5)}_{2}$}
\put(250,-60){$s^{(5)}_{3}$}
\put(365,-60){$s^{(5)}_{4}$}
\put(435,-60){$s^{(5)}_{5}$}
\put(345,7){$s^{(3)}_{3}s^{(3)}_{1}$}
\put(400,7){$s^{(3)}_{3}s^{(3)}_{2}$}
\put(453,7){$s^{(3)}_{3}s^{(3)}_{3}$}
}
%
%

\noindent
{\small
\setlength{\unitlength}{.22mm}
\begin{picture}(1000,300)(-25,-100)
\thicklines
\put(-30,170){\begin{minipage}[t]{2in}
%
%
\begin{fig}
\label{fig:btwo}
\end{fig}
\end{minipage}
}
\put(0,0){\mainfig}
\end{picture}
}

\noindent
By using a more rough analogy,
Figure \ref{fig:btwo} shows that an embedding is represented as a 
refinement of a partition
of a unit interval in the real line ${\bf R}$.
Then relations of inverse system 
in Theorem \ref{Thm:inverseb}(i)
mean that
the set of such refinements is a directed set.

\ssfr{Acknowledgment}
The author would like to Akira Asada for 
his idea of embeddings among Cuntz algebras.

%
%

%
\end{document}